\documentclass
{article}
\usepackage[T2A]{fontenc}
\usepackage[cp1251]{inputenc}
\usepackage[english,russian]{babel}
\usepackage[tbtags]{amsmath}
\usepackage{amsfonts,amssymb,mathrsfs,amscd}
\usepackage{wasysym}
\usepackage{verbatim}
\usepackage{eucal}
\usepackage{upgreek}
\usepackage{graphicx}
\usepackage{enumerate}
\let\No=\textnumero

\usepackage[hyper]{msb-a}
\JournalName{}

\theoremstyle{plain}

\newtheorem{propos}{Предложение}

\theoremstyle{definition}
\newtheorem{proof}{Доказательство}
\newtheorem{remark}{Замечание}
\newtheorem{definition}{Определение}

\renewcommand{\leq}{\leqslant} 
\renewcommand{\geq}{\geqslant}
\newcommand{\RR}{\mathbb{R}} 
\newcommand{\CC}{\mathbb{C}} 
\newcommand{\NN}{\mathbb{N}}

\DeclareMathOperator{\dist}{{\sf dist}}
\DeclareMathOperator{\har}{{\sf har}} 
\DeclareMathOperator{\conv}{{\sf conv}}

\begin{document} 
\title{О расстоянии Гарнака}
	
\author[B.\,N.~Khabibullin]{Б.\,Н.~Хабибуллин}
\address{Башкирский государственный университет, ИМВЦ УФИЦ РАН}
\email{khabib-bulat@mail.ru}

\udk{517.572}

 \maketitle

\begin{fulltext}

\begin{abstract} В заметке рассматриваются два типа оценок расстояния Гарнака в ограниченных  областях конечномерного евклидова пространства. Первый тип основан на геометрическом  понятии энтропия линейной связности.
Мы применяли   это понятие ранее для  исследования нулевых множеств голоморфных функций в весовых классах голоморфных функций. Второй тип основан на понятии показателя разделённости подмножества в области, которое вводится в этой заметке.  
Библиография:  10 названий 

Ключевые слова: гармоническая функция, неравенство Гарнака, расстояние Гарнака, линейная связность, разделённость множества в области





\end{abstract}

\markright{О расстоянии Гарнака}

	
$\NN:=\{1,2, \dots\}$, $\RR$, $\CC$ ---  множества {\it натуральных, вещественных, комплексных  чисел\/}
вместе с их стандартными порядковыми, алгебраическими, геометрическими, топологическими структурами;
$\RR^+:=\bigl\{x\in \RR\bigm| x\geq 0\bigr\}$ --- множество  {\it положительных чисел.\/}
Всюду далее ${\tt d}\in \NN$ --- размерность {\it евклидова пространства\/} 
$\RR^{\tt d}$ с {\it евклидовой нормой\/} $|x|:=\sqrt{x_1^2+\dots +x_{\tt d}^2}$ 
для $x:=(x_1,\dots ,x_{\tt d})\in \RR^{\tt d}$ и функцией $\dist(\cdot, \cdot)$  {\it евклидова расстояния}.
Пространство $\RR^{\tt 2}$, когда это уместно,  отождествляется с комплексной плоскостью $\CC\ni x+iy{\underset{x,y\in \RR}{\longmapsto}}(x,y)\in \RR^{\tt 2}$, если удаётся отвлечься от её комплексной структуры.

Через $\overline\RR:=\{-\infty\}\cup \RR\cup \{+\infty\}$ обозначаем {\it расширенную числовую прямую;\/}
а $\overline \RR^+:=\RR\cup \{+\infty \}$.
 При $r\in \overline\RR$ и  $x\in \RR^{\tt d}$ через $B_x(r)$ и  $\overline  B_x(r)$
обозначаем соответственно {\it открытый\/} и  {\it замкнутый шар радиуса $r$ с центром\/} $x$. 
Так, если $r\leq 0$, то  $B_x(r)$ --- {\it пустое множество\/} $\emptyset$, так же, как $\overline B_x(r)=\emptyset$  при   $r<0$. 
Всюду далее {\it положительность\/} --- это $\geq 0$, а {\it отрицательность\/} --- это $\leq 0$. 

Для подмножества $S\subset \RR^{\tt d}$
через $\overline S$,  $\complement S:=\RR^{\tt d}\setminus S$,
 и  $\partial S$ обозначаем соответственно {\it замыкание, дополнение\/} и {\it границу\/} множества $S$ в  
$\RR^{\tt d}$. Всюду   $D\neq \emptyset$ --- 
{\it область,\/} т.е. открытое связное подмножество  в $\RR^{\tt d}$. 
{\it Выпуклую оболочку\/} подмножества $S\subset \RR^{\tt d}$ обозначаем  через $\conv S$.
Таким образом, для  $x\in \RR^{\tt d}$ и $y\in \RR^{\tt d}$  $[x,y]:=\conv\{x,y\}$ --- {\it отрезок с концами $x$ и $y$.\/}
{\it Сегментальная  оболочка\/}  $\boxtimes S:=\cup \bigl\{ [x,y] \bigm| x,y \in S\bigr\}$
 --- это объединение всех {\it отрезков\/} $[x,y]$ {\it с концами\/} в точках $x$ и $y$ из множества $S$ 
 \cite[\S~2, 7]{BF},  \cite[3.1]{Aza79}, \cite[\S~11]{KhaKhaChe08I}, \cite{KhaKhaChe08II}. 
{\it Звёздная  оболочка множества $S$  с началом в точке $z\in D$} --- это объединение 
 $\hexstar_{z} S:=\bigcup \bigl\{ [z,x] \bigm| x \in S \bigr\} $  всех отрезков, соединяющих  $z$ со всеми точками $x\in S$. 
Через $\diameter\! S:=\sup\limits_{x,y\in S}|x-y|$ обозначаем
{\it евклидов диаметр} множества $S$.

  $\har(D)$ ---  векторное пространство над $\RR$ {\it гармонических функций\/} на $D$, а 
$\har^+(D)\subset \har(D)$ --- конус {\it положительных\/} гармонических функций.

\begin{definition}[{\cite[1.3]{Rans}, \cite{Kohn}}]
 Для  области $D\subset \RR^{\tt d}$ {\it расстояние  Гарнака\/} между точками $x\in D$ и $y\in D$ 
на  $D$  обозначаем через  
\begin{equation}\label{hard}
\hspace{-1mm}\dist_{\har}^D (x,y):=\inf\Bigl\{ d\in \RR^+\Bigm|
\frac{1}{d}  h(y)\leq h(x) \leq d h(y)\quad  
\forall h\in \har^+(D)\Bigr\}.
\end{equation}
\end{definition}
Основные свойства расстояния Гарнака \eqref{hard} отражены в   \cite[1.3]{Rans}, \cite{Kohn}. 
  Например,  $\ln \dist_{\har}^D$ --- {\it полуметрика на $D$, непрерывная  относительно евклидовой метрики в $D$
\cite[\S~1]{Kohn},} которая совпадает с  {\it гиперболической метрикой\/} в случае  {\it односвязной\/} области 
 $D\subset \CC$, отличной от $\CC$ \cite{Chi12}. 
Наряду  с {\it симметричностью\/} имеет место {\it мультипликативное  неравенство треугольника\/}
\begin{equation}\label{pmt}
\dist_{\har}^D(x,y)\leq \dist_{\har}^D(x,z)\cdot \dist_{\har}^D(z,y)\text{ \it для любых $x,y,z\in D$}.
\end{equation}  
Для расстояния Гарнака справедлив  принцип подчинения  \cite[теорема 1.3.6]{Rans}, \cite[\S~3, 
теорема 3.3]{Kohn}, который будет использован в следующей упрощённой  форме
\cite[следствие  1.3.7]{Rans}, \cite[\S~3, 
теорема 3.2]{Kohn}: 
\begin{equation}\label{pr}
\dist_{\har}^D\leq \dist_{\har}^{D'} \quad \text{\it  на $D'$ для областей  $D'\subset D$}.
\end{equation} 
Для шара расстояние Гарнака выписывается явно \cite[\S~3, приложение]{Kohn}:
\begin{equation}\label{DB}
\dist_{\har}^{B_{x_0}(R)}(x_0, x)=\frac{\bigl(R+|x-x_0|\bigr)R^{{\tt d-2}}}{\bigl(R-|x-x_0|\bigr)^{{\tt d}-1}}
\text{ \it при  $x\in B_{x_0}(R)$},
\end{equation}
где правая часть возрастает при увеличении расстояния $|x-x_0|$ от точки $x$ до центра $x_0$ шара $B_{x_0}(R)$.

\section{Расстояние Гарнака и энтропия линейной связности}\label{eac} 

\begin{definition}[{\cite[определение 2.1]{KhaKhaChe08I}, \cite[определение III.1]{KhaKhaChe08II}}]\label{elc}
{\it Энтропией линейной связности\/} подмножества $S$ в  области $D\subset \RR^{\tt d}$
называем величину
\begin{equation}\label{els}
{\sf eac}_D(S):=\sup_{x,y\in S} \inf_{l(x,y)\subset D}\frac{\bigl|l(x,y)\bigr|}{\dist\bigl(l(x,y), \complement D\bigr)}\in \overline \RR^+,
\end{equation}
где  $l(x,y)\subset D$ --- всевозможные  {\it спрямляемые кривые  с концами $x$ и $y$ длины  
$\bigl|l(x,y)\bigr|$  в евклидовой метрике.\/} Спрямляемые кривые  с концами $x$ и $y$  здесь можно заменить на 
{\it ломаные, соединяющие $x$ и $y$}.
\end{definition}
 Случай ограниченной области $D\subset \CC$ детально рассмотрен и существенно используется  в  \cite[\S~2]{KhaKhaChe08I},  \cite[гл.~III]{KhaKhaChe08II}. Отметим дополнительные элементарные  свойства энтропии линейной связности, которые приведены  в \cite[предложение 2.3]{KhaKhaChe08I} для областей в $\CC$, но без изменений переносятся и на области в $\RR^{\tt d}$ \cite[\S~2, замечание]{KhaKhaChe08I}.  
Далее всюду в \S~\ref{eac} области $D$  и $D'${\it ограниченные.} 

\begin{propos} ${\sf eac}_D (S) =  {\sf eac}_D (D\cap \overline S)$ и имеют место два принципа подчинения: 
 eсли  $S \subset S'\subset D$, то ${\sf eac}_D (S)\leq {\sf eac}_D (S')$, а  если $D\subset D'\subset \mathbb R^{\tt d}$,
то ${\sf eac}_D (S) \geq {\sf eac}_{D'} (S)$. Кроме того, ${\sf eac}_D (S) < +\infty$, если и только если  
$\overline S\subset D$.
\end{propos}

Доказательство следующего утверждения также  без изменений переносится на $\RR^{\tt d}$  из доказательства \cite[предложение 2.2]{KhaKhaChe08I} для $\CC$. 

\begin{propos}\label{prop:mot} Если ${\sf eac}_D (S)<C < +\infty$, то для любой пары  $x,y\in S$ найдется конечная последовательность шаров  $B_{x_k}(r)\subset D$, $k= 0, \dots , l+1$,
 для которых одновременно выполнены два условия:
\begin{enumerate}[{\rm (Bi)}] 
\item\label{D1} $x=x_0$, $y=x_{l+1}$ и $|x_{k-1}-x_k|\leq r/2$ при всех $k=1, \dots , l+1$;
\item\label{D2} имеет место оценка $l\leq 2C$.
\end{enumerate}
\end{propos}

\begin{propos}\label{pr:adH}
Для  подмножества $S$  ограниченной области $D\subset \RR^{\tt d}$
\begin{equation}\label{end:aSD}
\sup_{x,y \in S} \dist_{\har}^D(x, y)\leq (3\cdot 2^{\tt d-2})^{2{\sf eac}_D (S)+1}
\leq 2^{2{\tt d}({\sf eac}_D (S)+1)}.
\end{equation}
\end{propos}
Доказательство также  переносится на $\RR^{\tt d}$  из \cite[предложение 3.2]{KhaKhaChe08I}. Но мы приведём его полностью для контроля над числовыми параметрами.
\begin{proof}  Можно считать энтропию ${\sf eac}_D (S)$ конечной.

Пусть $x, y \in S$ и ${\sf eac}_D (S)<C < +\infty$. Рассмотрим конечную
последовательность кругов $B_{x_k}(r)\subset D$, $k=0, \dots , l+1$,
из предложения  \ref{prop:mot}, удовлетворяющих условиям (B\ref{D1})--(B\ref{D2}).
Из свойства \eqref{pmt} расстояния Гарнака следует
$$
\dist_{\har}^D(x, y)\leq \prod_{k=1}^{n+1} \dist_{\har}^D(x_{k-1}, x_k).
$$
Отсюда, применяя к правой части  принцип подчинения \eqref{pr} при
$B_{x_k}(r)\subset D$ и каждом   $k=0, \dots , l+1$, получаем
\begin{multline*}
\dist_{\har}^D(x, y)\leq \prod_{k=1}^{l+1} \dist_{\har}^D(x_{k-1}, x_k)\leq
\prod_{k=1}^{l+1} \dist_{\har}^{B_{x_k}(r)}(x_{k-1}, x_k).
\\
\overset{\eqref{DB}}{=} 
\prod_{k=1}^{l+1}\frac{\bigl(r+|x_{k-1}-x_k|\bigr)r^{{\tt d-2}}}{\bigl(r-|x_{k-1}-x_k|\bigr)^{{\tt d}-1}}
\leq \prod_{k=1}^{l+1}\frac{(r+r/2)r^{{\tt d-2}}}{(r-r/2)^{{\tt d}-1}}
=\prod_{k=1}^{l+1}3\cdot 2^{\tt d-2}=(3\cdot 2^{\tt d-2})^{l+1}
\end{multline*}
Для получения из последней оценки  соотношения \eqref{end:aSD} остается
обратиться к  оценке $l\leq 2C$ из (B\ref{D2}) и учесть произвол в выборе
$C>{\sf eac}_D (S)$.
\end{proof}

Для специальных подмножеств в областях 
энтропия линейной связности часто может оцениваться достаточно просто.



\begin{propos}\label{eq:diamSS}
Если выпуклая оболочка $\conv S$ множества $S$ содержится в ограниченной области $D\subset \RR^{\tt d}$, то 
\begin{equation}\label{S1} 
{\sf eac}_D (S) \leq \dfrac{\diameter\! S}{\dist(\conv S, \complement D)}.
\end{equation}
Если  сегментальная  оболочка $\boxtimes S$
 содержится в $D$, то 
\begin{equation}\label{S2} 
{\sf eac}_D (S) \leq \dfrac{\diameter\! S}{\dist(\boxtimes  S, \complement D)}.
\end{equation}
Если звёздная  оболочка  $\hexstar_{z} S$  с началом $z$ содержится в $D$,   то 
\begin{equation}\label{S3} 
{\sf eac}_D (S) \leq \dfrac{\diameter\! S}{\dist(\hexstar_{z} S, \complement D)}.
\end{equation}
\end{propos}
\begin{proof} Для оценки сверху энтропии линейной связности в качестве кривых $l(x,y)$ в определении 
\eqref{els} в первых двух  случаях   достаточно рассмотреть отрезки, соединяющие точки $x$ и  $y$, а 
в третьем случае --- объединения пар отрезков, соединяющих $x$ с $z$ и $z$ с $y$. 
\end{proof}

\begin{remark} В \cite[\S\S~9--11, предложение 11.3]{KhaKhaChe08I} имеются и другие более тонкие результаты об  энтропии линейной связности для ограниченных областей $D\subset \CC$, которые могут  быть использованы для оценок  расстояния Гарнака.  
\end{remark}

\section{Расстояние Гарнака и показатель разделённости}
По определению \eqref{els} энтропия линейной связности  --- это одно значение, или параметр. Замена кривых $l(x,y)\subset D$ ломаными $l(x,y)$ с конечным числом вершин  в определении  
\eqref{els} не меняет значения энтропией линейной связности ${\sf eac}_D (S)$.
Если  допустить большее число  параметров, характеризующих подмножества $S$ в области $D$, то
 можно расширить  оценки  расстояния Гарнака. Так, можно развить содержание предложения 
\ref{prop:mot} и, как следствие, предложений \ref{pr:adH}--\ref{eq:diamSS}, допуская переменные радиусы шаров $B_{x_k}( \cdot )$, конечные цепи которых соединяют точки $x_0\in  D$ и $x\in S\subset D$. Именно таким путём определяется метрика Кабояши, напрямую связанная с расстоянием Гарнака, во всяком случае, в $\CC$ \cite{Chi12}. 
 Другой путь ---  ограничиться условиями на последовательность их центров $x_k$, как в (B\ref{D1}), отвлекаясь от шаров. Ниже в предложениях  \ref{prB}--\ref{lpr} 
реализуется вариант таких  двухпараметрических  условий.

\begin{propos}\label{prB} 
Пусть   область  $D\subset \RR^{\tt d}$ отлична от $\RR^{\tt d}$, 
а для пары точек $x\in D$ и $y\in D$
отношение расстояния между ними  к сумме расстояний от них до границы области $D$
не больше $q\in [0,1)$:
\begin{equation}\label{otx+}
\frac{|x-y|}{\dist(x, \complement D)+\dist(y, \complement D)}\leq q<1.
\end{equation}
 Тогда 
\begin{equation}\label{z}
\dist_{\har}^D(x,y)\leq \frac{2^{2{\tt d}}}{(1-q)^{2({\tt d-1})}}.
\end{equation}
\end{propos}

\begin{proof} 
Положим 
\begin{equation}\label{otx}
\begin{cases}
r:=\dist(x, \complement D)>0,\\
R:=\dist(y, \complement D)>0,
\end{cases}
\quad \overset{\eqref{otx+}}{\Longrightarrow} \quad  \frac{|x-y|}{r+R}\leq q<1.
\end{equation}
Согласно \eqref{otx}   шары $B_x(r)\subset D$ и $B_y(R)\subset D$ имеют непустое пересечение, а их объединение и, тем более, отрезок $[x,y]$ содержаться  в $D$. Кроме того, по определению $r$ и $R$ в \eqref{otx},  замыкание ни одного из этой пары шаров не может содержаться в другом шаре, откуда  $|x-y|\geq |R-r|$. 
Опишем выбор специальной точки $z$, лежащей на отрезке $[x,y]$ и принадлежащей обоим шарам.
Для этого отложим от концов $x$ и $y$ отрезка $[x,y]$ два перекрывающихся отрезка 
$[x,y_1]$ и $[x_1, y]$ на прямой, проходящей через $x$ и $y$, длины соответственно  $|y_1-x|=r$ и $|y-x_1|=R$ . Тогда 
\begin{multline*}
\min \biggl\{\frac{\bigl|[x_1,y_1]\bigr|}{r}, \frac{\bigl|[x_1,y_1]\bigr|}{R}\biggr\}\geq \frac{\bigl|[x_1,y_1]\bigr|}{R+r}=\frac{|y_1-x_1|}{R+r}\\
=
\frac{(R+r)-|x-y|}{R+r}=1-\frac{|x-y|}{R+r}\overset{\eqref{otx}}{\geq} 1-q. 
\end{multline*}
Выбираем  точку $z:=\frac12x_1+\frac12y_1$ как середину отрезка $[x_1,y_1]$. По построению и ввиду неравенства 
$|x-y|\geq |R-r|$ такая точка принадлежит отрезку $[x,y]$, каждому из шаров и  одновременно отрезкам  $[x,y_1]$
и $[x_1,y]$, а следовательно,
\begin{multline*}
\max\biggl\{\frac{\bigl|[x,z]\bigr|}{r}, \frac{\bigl|[z,y]\bigr|}{R}\biggr\}
=\max\biggl\{\frac{r-\bigl|[z,y_1]\bigr|}{r}, \frac{R-\bigl|[x_1, z]\bigr|}{R}\biggr\}
\\=
1-\frac12 \min \biggl\{\frac{\bigl|[x_1,y_1]\bigr|}{r}, \frac{\bigl|[x_1,y_1]\bigr|}{R}\biggr\}
\leq 1-\frac12(1-q)=\frac{1+q}{2}<1.
\end{multline*}
Отсюда по принципу подчинения \eqref{pr} и явному виду \eqref{DB} расстояния Гарнака 
для шара $B_x(r)$, содержащего точку $z$, получаем
\begin{multline*}
\dist_{\har}^D(x,z)\leq \dist_{\har}^{B_x(r)}(x,z)
\overset{\eqref{DB}}{=}
\frac{\Bigl(r+\bigl|[x,z]\bigr|\Bigr)r^{{\tt d-2}}}{\Bigl(r-\bigl|[x,z]\bigr|\Bigr)^{\tt d-1}}
\\
=\frac{1+\bigl|[x,z]\bigr|/r}{\Bigl(1-\bigl|[x,z]\bigr|/r\Bigr)^{\tt d-1}}
\leq \frac{1+(1+q)/2}{\bigl(1-(1+q)/2\bigr)^{\tt d-1}}
=2^{\tt d-2}\frac{3+q}{(1-q)^{\tt d-1}}, 
\end{multline*}
а также, из очевидных соображений симметрии, такую же оценку для 
$$
\dist_{\har}^D(z,y)\leq 2^{\tt d-2}\frac{3+q}{(1-q)^{\tt d-1}}.
$$
Исходя из этих неравенств, по  мультипликативному правилу треугольника \eqref{pmt} 
\begin{multline*}
\dist_{\har}^D(x,y)\overset{\eqref{pmt}}{\leq} \dist_{\har}^D(x,z)\cdot \dist_{\har}^D(z,y)
\leq \Bigl(2^{\tt d-2}\frac{3+q}{(1-q)^{\tt d-1}}\Bigr)^2\overset{\eqref{otx+}}{=} \frac{2^{2{\tt d}}}{(1-q)^{2({\tt d-1})}}
\end{multline*}
 что, ввиду ограничений $0\leq q\overset{\eqref{otx}}{<}1$, даёт требуемую оценку  \eqref{z}.
\end{proof}

Чем больше   левая часть \eqref{otx+}, тем больше отделены друг от друга точки $x\in D$
и $y\in D$ внутри  $D$ относительно $\complement D$ или $\partial D$. Это мотивирует 

\begin{definition}\label{betw} Пусть $D$ --- область в $\RR^{\tt d}$.
{\it Разделённостью ${\sf sep}_D(x,y)\in \RR^+$ точек $x\in D$ и $y\in D$ внутри}  $D$ называем  дробь 
из левой  части 
\eqref{otx+}:
\begin{equation}\label{otx+0}
{\sf sep}_D(x,y):=\frac{|x-y|}{\dist(x, \complement D)+\dist(y, \complement D)}.
\end{equation}
{\it Разделённость  ${\sf sep}_D(x_0,\dots, x_l)$ последовательности точек  $x_0, x_1, \dots , x_l\in D$
внутри $D$\/} ---  максимум разделённостей  по всем парам  соседних точек:
\begin{equation}\label{bwX}
{\sf sep}_D(x_0,\dots, x_l):=\max_{k=1, \dots, l}
{\sf sep}_D(x_{k-1}, x_k)\in  \RR^+.
\end{equation}
{\it Разделённость ${\sf sep}_D(S; x_0,l)$ подмножества $S\subset D$ при начале $x_0\in D$ по длине $l\in \NN$ внутри $D$} 
определяем как 
\begin{equation}\label{DL}
{\sf sep}_D(S; x_0,l)
:=\sup_{x_l\in S}\, \inf_{x_1, \dots , x_{l-1}\in D} \, {\sf sep}_D(x_0, \dots, x_l)\in \overline \RR^+.
\end{equation}
\end{definition}
Из определения \ref{betw} сразу следует 
\begin{propos}\label{pBB}
Если ${\sf sep}_D(x_0,\dots, x_l)\overset{\eqref{bwX}}{\leq} 1$, то ломаная 
$$\bigcup_{k=1}^l[x_{k-1}, x_k] $$
с  последовательными  вершинами $x_k$, $k=0, 1, \dots, l$, и даже цепочка шаров 
\begin{equation}\label{lom}
B_{x_k}(r_k), \quad r_k:=\dist(x_k, \complement D), \quad k=0, 1, \dots, l
\end{equation}
лежат в $D$, а  ${\sf sep}_D(x,x)=0$ для любой точки $x\in D$.
\end{propos} 
 
По предложению  \ref{pBB} при заданном числе $q<1$  ключевое  далее  условие 
\begin{equation}\label{DLs}
{\sf sep}_D(S; x_0,l)\leq q<1
\end{equation}
можно записать в двух эквивалентных формах  без упоминания о понятии 
разделённости  множества $S\subset D$ при начале $x_0\in D$ по длине $l\in \NN$ внутри $D$:

\begin{itemize}
\item[{$[\between_q^l]$}] {\it Для  каждых $x\in S$  и $q'>q$ найдётся 
ломаная в $D$ с концами  $x_0$ и $x$  из $l$ звеньев с длиной каждого звена не больше домноженной на $q'$  суммы двух расстояний от концов этого звена до границы $\partial D$.}
\item[{$[\between_q^l]$}]
{\it Для каждых   $x\in S$  и $q'>q$ найдётся  цепь  из $l+1$ шаров} 
\begin{equation}\label{Hss}
B_{x_k}(r_k),  \quad x_k\in D, \quad  k=0, \dots , l,  
\end{equation}
{\it с последним центром $x_l=x$ и радиусами 
$r_k=\dist(x_k, \partial D)$, для которых\/ 
$|x_k-x_{k+1}|\leq q'(r_k+r_{k+1})$ при каждом\/} $k=0, 1, \dots, l-1$.
\end{itemize}

\begin{propos}\label{lpr}
Пусть  $D$ --- область в $\RR^{\tt d}\neq D$, $x_0\in D$, $l\in \NN$,  а  подмножество $S\subset D$ удовлетворяет 
условию \eqref{DLs}, т.е. обладает  свойством   $[\between_q^l]$  с $q<1$ в любой из двух эквивалентных форм. Тогда   
 \begin{equation}\label{ita}
\sup_{x\in S} \dist_{\har}^D(x_0,x) \leq \frac{2^{2{\tt d}l}}{(1-q)^{({\tt d-1})l}}. 
\end{equation}
\end{propos}
\begin{proof} Применяем   предложение \ref{prB} к каждой паре последовательных вершин ломаной из свойства 
$[\between_q^l]$, затем не более чем  $l$ раз мультипликативное правило треугольника \eqref{pmt}, после чего в полученном неравенстве  переходим  к точной верхней грани по $x_l\in S$, как в  определении \eqref{DL}, и, наконец, учитываем произвол в выборе числа $q'>q$ для $q<1$. 
\end{proof}

\begin{remark} 
Для свойства  $[\between_q^l]$ во второй трактовке  в терминах шаров  прослеживаются явные 
параллели  с так называемыми цепями Гарнака из шаров \eqref{Hss}, использованными в \cite[3]{JerKen82}, 
 \cite[гл.~3, 1]{CapKenLan05} при  определении некасательно достижимых (NTA-)областей для  исследования гармонической меры и граничных свойств функций.  Наряду с одновременно и внутренним, и внешним требованием к границе NTA-области $D$, известным как  условие штопора  
  --- ``corkscrew condition'' 
 \cite[(3.1)--(3.2)]{JerKen82},  \cite[определение 3.1, (i)--(ii)]{CapKenLan05}, \cite[определение 1.5]{HWcorkscrew},  требовалось  и существование цепей последовательно пересекающихся   $l+1$ шаров \eqref{Hss}, соединяющих 
точки  $P_1\in B_{x_0}(r_0)$ и $P_2\in B_{x_{l}}(r_{l})$, в некотором смысле $(\varepsilon,R_0)$-близкие к граничным точкам области  $D$  \cite[(3.3)]{JerKen82}, \cite[определение 3.1, (iii)]{CapKenLan05}, но с иными условиями на радиусы 
и центры шаров, которые  должны быть $M$-некасательными для 
числа $M>1$ в том смысле, что 
\begin{equation}\label{cB+}
\frac1{M} r_k\leq \dist\bigl(B_{x_k}(r_k), \partial D\bigr)\leq Mr_k
\quad\text{\it для всех $k=0, \dots, l$}, 
\end{equation}
а также с ограничением на диаметры шаров 
\begin{equation}\label{cB+1}
\min_{k}\diameter\! B_{x_k}(r_k)\geq  \frac1{M}\min \{\dist (x_0, \partial D), \dist (x_{l},\partial D)\} .
\end{equation}
Такого типа  условия  тоже применимы к оценкам расстояния Гарнака, но со значительным  усложнением правой части \eqref{ita} дополнительными множителями  и слагаемыми, выражаемыми через 
$M$ и $R_0$,  к тому же    на  жёстком фоне условий  штопора для границы $\partial D$ и требований  \eqref{cB+} и \eqref{cB+1} соизмеримости центров и  радиусов шаров с близостью к границе $\partial D$. 
\end{remark}
\end{fulltext}

\end{document}